\input ./style/arxiv-general.cfg
\documentclass[aap,MSNbibl,seceqn,dvips]{arximspdf}
\makeatletter
   \@ifpackageloaded{graphicx}{}{\usepackage{graphicx}}
\makeatother
\usepackage{algorithm}
\usepackage{algpseudocode}

\doi{10.1214/14-AAP1042} %kopijuoti is PTS
\volume{25}
\issue{4}
\pubyear{2015}
\firstpage{2039}
\lastpage{2065}
\docsubty{FLA}

\makeatletter
\newcommand{\dd}{\mathrm{d}}
\newcommand{\rrvert}{\vert}
\newcommand{\llvert}{\vert}
\newcommand{\eqref}[1]{(\ref{#1})}
\newproclaim{defn}{Definition}
\newtheorem{teo}{Theorem}
\newproclaim{rem}{Remark}
\newtheorem{cor}{Corollary}
\newtheorem{lem}{Lemma}
\newproclaim{assumption}{Assumption}
\newtheorem{prop}{Proposition}
\newproclaim{example}{Example}
\makeatother

\begin{document}
\begin{frontmatter}

\title{An integral equation for Root's barrier and the generation of
Brownian increments}
\runtitle{An integral equation for Root's barrier}

\begin{aug}
\author[A]{\fnms{Paul}~\snm{Gassiat}\corref{}\ead[label=e1]{gassiat@math.tu-berlin.de}\thanksref{T1}},
\author[B]{\fnms{Aleksandar}~\snm{Mijatovi\'{c}}\ead[label=e2]{a.mijatovic@imperial.ac.uk}\thanksref{T2}}\\
\and
\author[C]{\fnms{Harald}~\snm{Oberhauser}\ead[label=e3]{harald.oberhauser@oxford-man.ox.ac.uk}\thanksref{T1,T3}}
\runauthor{P. Gassiat, A. Mijatovi\'{c} and H. Oberhauser}
\affiliation{Technische Universit\"at Berlin, Imperial College and
University of Oxford}
%\dedicated{}
\address[A]{P. Gassiat\\
Institut f\"{u}r Mathematik\\
Technische Universit\"at Berlin\\
Strasse des 17.~Juni 136\\
10623 Berlin\\
Germany\\
\printead{e1}} %adresu isvedimo komanda gale!
\address[B]{A. Mijatovi\'{c}\\
Department of Mathematics\\
Imperial College\\
180 Queen's Gate\\
London SW7 2AZ\\
United Kingdom\\
\printead{e2}}
\address[C]{H. Oberhauser\\
Oxford--Man Institute\\
University of Oxford\\
Walton Well Road\\
Oxford OX2 6ED\\
United Kingdom\\
\printead{e3}}
\end{aug}
\thankstext{T1}{Supported in part by the European Research Council, Grant agreement Nr.~258237.}
\thankstext{T2}{Supported in part by the Humboldt Foundation Research Fellowship, GRO/1151787 STP.}
\thankstext{T3}{Supported in part by the European Research Council, Grant agreement Nr.~291244.}

% HISTORY:
\received{\smonth{10} \syear{2013}}
\revised{\smonth{4} \syear{2014}}
%\accepted{\smonth{} \syear{}}

% ABSTRACT
%
\begin{abstract}
We derive a nonlinear integral equation to calculate Root's solution
of the Skorokhod embedding problem for atom-free target measures.
We then use this to efficiently generate bounded time--space increments
of Brownian motion and give a parabolic version of Muller's classic
``Random walk over spheres'' algorithm.
\end{abstract}

% KEYWORDS
\end{frontmatter}

\setcounter{footnote}{3}
%s1 #&#
\section{Introduction}\label{sec1}

Let $\mu$ be a zero-mean probability measure on the real line and
$B= (B_{t} ){}_{t\geq0}$ denote a one-dimensional Brownian
motion. The Skorokhod embedding problem given by $\mu$ consists of
constructing a stopping time $\tau$ such that
{\renewcommand{\theequation}{SEP$_\mu$}{
%
%e1.1 #&#
\begin{equation}
\label{SEP} \qquad B_{\tau}\sim\mu\quad\mbox{and}\quad B^{\tau}=
(B_{t\wedge\tau} )_{t\geq0}\qquad\mbox{is uniformly integrable.}
\end{equation}}}\setcounter{equation}{0}%
More than 50 years after Skorokhod \cite{MR0185619}, we can now choose
from a wide range of different stopping times which solve this
problem~\cite{Hobson2011,obloj2004skorokhod}.
In general such a stopping time may depend in a very complex way on
the Brownian trajectory. This can make it computationally expensive
(or even intractable) in applications to determine the actual realisation
of the stopping time $\tau$ for a given Brownian trajectory. From
this point of view, one of the earliest solutions of~(\ref{SEP}),
the so-called \textit{Root solution}, is one that stands out: in 1969
Root~\cite{MR0238394} showed that if $\mu$ has zero-mean and a second
moment, then there exists a closed subset of time--space,
the so-called \textit{Root barrier},
\[
R\subset[0,\infty]\times[-\infty,\infty],
\]
such that the hitting time
\[
\tau=\inf\bigl\{ t>0\dvtx  (t,B_{t} )\in R \bigr\} \qquad(\inf\varnothing=
\infty)
\]
solves (\ref{SEP}) given by $\mu$. The Root barrier $R$ can be
described by a lower semicontinuous barrier function $r$,
\[
R= \bigl\{ (t,x )\dvtx t\geq r (x ) \bigr\},
\]
and, among all solutions $\tilde{\tau}$ of~(\ref{SEP}), it has
the key property of minimising $\mathbb{E} [\tilde{\tau}^{2} ]$; see
Rost~\cite{MR0445600} and Loynes~\cite{MR0292170}. Unfortunately,
Root's existence proof is not constructive, and until recently it was
not known how to characterise or compute $R$ (or, equivalently,~$r$)
in terms of the measure $\mu$. A seminal paper by Hobson \cite{hobson1998robust}
on applications to model independent hedging of exotic options led
to a revived interest in~(\ref{SEP}) (the Root solution gives a
lower bound on options on variance), and motivated by such applications,
the Root barrier was more recently identified as the free boundary
of a parabolic obstacle problem (work of Dupire, Cox and Wang,
Oberhauser and Reis, \cite
{Dupire2005,cox2011root,oberhauserdosreis12,2013arXiv13084363C}).
This allows one to compute $R$ in two steps: firstly, solve numerically
the nonlinear PDE (using finite difference or BSDE methods), and secondly,
numerically calculate the associated free boundary of this~PDE.

The first and main contribution of this paper consists of characterising
the barrier function $r$ directly via a nonlinear integral equation.
More precisely, if $\mu$ is atom-free, then $r$ solves the following
equation:
%
%e1.1 #&#
%e1.2 #&#
\begin{eqnarray}\label{eqintegralequ}
\qquad u_{\delta_{0}} (x )-u_{\mu} (x )=g \bigl(r (x ),x \bigr)-\int
_{ \{ y\dvtx r (y )<r (x ) \} }g \bigl(r (x )-r (y ),x-y \bigr)\mu(\dd y )
\nonumber\\[-8pt]\\[-12pt]
\eqntext{\forall x\in(-\infty,\infty).}
\end{eqnarray}
Here $g (t,x )=\sqrt{\frac{2t}{\pi}}e^{-{x^{2}}/{(2t)}}-\llvert
x\rrvert\operatorname{Erfc} (\frac{\llvert x\rrvert}{\sqrt{2t}}
)=\mathbb{E}L_{t}^{x}$
where $(L_t^x)_{t,x}$ is the Brownian local time, and $u_{\mu}$
and $u_{\delta_{0}}$ are the potential functions%
\footnote{That is, $u_\mu(x) = -\int|y-x| \mu(\dd y)$ is the formal
density of the occupation measure $\mu U=\int_{0}^{\infty}\mu
P_{t}\,\dd t$
%and $\delta_{0}U=\int_{0}^{\infty}\delta_{0}P_{t}\dd t$
where $P_{t}$ denotes the semigroup of Brownian motion.%
}
of the measures $\mu$ and the Dirac delta $\delta_{0}$, respectively.
The derivation of this integral equation is short, intuitive, and entirely
probabilistic as it relies solely on the It\^{o}--Tanaka formula and the
fact that the local time is an additive functional of the path of
Brownian motion.

%Note that all terms in equation~(\ref{eqintegralequ}) can be explicitly
%expressed as functions of $x,y,r\left(x\right),r\left(y\right)$,
%making~(\ref{eqintegralequ}) a nonlinear integral equation for
%the function $r$, with the integral kernel given by $g\left(t,x\right)=
%\mathbb{E}\left[L_{t}^{x}\right]$
%(see Corollary~\ref{corintequation} for the explicit formula for
%$g$).
It is well known (see, e.g., \cite{peskir2005american}) that
the question of uniqueness of solutions of such nonlinear integral
equations is delicate
in general. In this case we give a short proof of the uniqueness of
the solution of~(\ref{eqintegralequ}) that applies to the class
of measures with a continuous barrier function via the uniqueness
of the viscosity solution of a nonlinear PDE characterising the Root
solution of~(\ref{SEP}) given in~\cite{oberhauserdosreis12}.

In the rest of the article we then specialise to the case of barriers
that have a barrier function that is symmetric around $0$, continuous,
and monotone. In this case it becomes numerically much easier to solve
(\ref{eqintegralequ}) since $r$ does not appear anymore in the
domain of the integral, and (\ref{eqintegralequ}) becomes a Volterra
type integral equation of the first kind. Furthermore, we again use
the viscosity approach of \cite{oberhauserdosreis12} to establish
sufficient and easy to verify conditions on the probability measure
$\mu$ which guarantee that its barrier has these properties. These
results give a theoretical justification for the application of a
simple numerical scheme to this integral equation, yielding a much
faster and more accurate numerical method for directly computing $r$,
for a class of symmetric probability measures $\mu$ with compact
support, than the nonlinear PDE approach described above.

The second contribution of this paper is to show that~(\ref{SEP}),
and in particular the Root solution described by the equation~(\ref
{eqintegralequ}),
can be very useful in sampling Brownian increments, an essential task
in Monte Carlo schemes. Recall the arguably simplest algorithm $ (\tau
_{0}^{\mathrm{sim}},X_{0}^{\mathrm{sim}} )= (0,0 )$
and
\[
\cases{ \displaystyle X_{n+1}^{\mathrm{sim}} = X_{n}^{\mathrm{sim}}+N_{n},
&\quad with i.i.d. $N_{n}\sim\mathcal{N} (0,1 )$,
\vspace*{3pt}\cr
\displaystyle
\tau_{n+1}^{\mathrm{sim}} = \tau_{n}^{\mathrm{sim}}+1.}
\]
Then the equality holds $ (\tau_{n}^{\mathrm{sim}},X_{n}^{\mathrm{sim}} )_{n\in\mathbb
{N}}\stackrel{\mathrm{Law}}{=} (\tau_{n},B_{\tau_{n}} )_{n\in\mathbb{N}}$
where $\tau_{n}=n$. We would like to stress here that this algorithm
works because $\tau_{1}$ solves~(\ref{SEP}): $B_{\tau_{1}}\sim\mathcal
{N} (0,1 )$.
In fact, setting $r\equiv1$, that is,
\[
R= \bigl\{ (t,x )\dvtx t\geq1,x\in[-\infty,\infty] \bigr\},
\]
it follows $\tau_{1}\equiv1=\inf\{ t>0\dvtx  (t,B_{t} )\in R \}$,
and we see that Root's solution for $\mu=\mathcal{N} (0,1 )$
yields the classical Euler scheme. Note, however, that at least in
principle, Root's result allows us to choose $\mu$ to be any probability
measure on real numbers. The canonical choice, as pointed out by Dupire,
in terms of speed of simulation on a standard computer, which is very
efficient in drawing quasi-random numbers from the uniform distribution,
is to take $\mu=\mathcal{U} [-1,1 ]$. In this case the barrier
function $r$ can be computed (once!) arbitrarily accurately via~(\ref
{eqintegralequ}),
yielding a simulation algorithm
\[
\cases{ \displaystyle X_{n+1}^{\mathrm{sim}} = X_{n}^{\mathrm{sim}}+U_{n},
&\quad with i.i.d. $U_{n}\sim\mathcal{U} [-1,1 ]$,
\vspace*{3pt}\cr
\displaystyle
\tau_{n+1}^{\mathrm{sim}} = \tau_{n}^{\mathrm{sim}}+r(U_{n}
).}
\]
Again we have $ (\tau_{n}^{\mathrm{sim}},X_{n}^{\mathrm{sim}} )_{n\in\mathbb{N}}\stackrel
{\mathrm{Law}}{=} (\tau_{n},B_{\tau_{n}} )_{n\in\mathbb{N}}$
where $ (\tau_{n} )_{n}$ denote the first hitting times
of $t\mapsto(t,B_{t} )$ of the Root barrier $R$, that is,~$\tau
_{1}=\inf\{ t>0\dvtx  (t,B_{t} )\in R \} $,
$\tau_{2}=\inf\{ t>\tau_{1}\dvtx  (t-\tau_{1},B_{t}-B_{\tau_{1}} )\in R \} $,
etc. What makes this algorithm particularly interesting, besides its
computational efficiency, is the fact that the time--space process
$ (t,B_{t} ){}_{t\geq0}$, and in particular the Brownian
motion itself, is uniformly bounded between consecutive sampling times
$\tau_{n}$ and $\tau_{n+1}$ for all $n\in\mathbb{N}$,
\[
\sup_{t\in[\tau_{n},\tau_{n+1} ]}\llvert B_{t}-B_{\tau_{n}}\rrvert<2
\quad\mbox{and}\quad\tau_{n+1}-\tau_{n}<\frac{2}{\pi}
\]
(the first inequality is sharp but the second is not; see Corollary~\ref
{corboundsuniform}).
Such a property is particularly useful in Monte Carlo schemes for
computing solutions of PDEs with time-dependent boundaries; similar
observations have been made by many different authors before, for
example, Milstein and Tretyakov,
Deaconu and Hermann, Deaconu, Lejay, and Zein \cite
{MR1722281,deaconeherrman2013,MR2673769},
by using different shapes (e.g.,~parallelepipeds) than $R$; however, the
above approach via the (\ref{SEP}) is extremal among these solutions
in the sense that it allows one to sample from the arguably simplest
distribution
for computational purposes $\mathcal{U} [-1,1 ]$. It is
also clear that Brownian scaling can be used to modify the above algorithm,
which is described in detail in Section~\ref{secgenerating-bounded-brownian},
to sample increments during which the uniform bound is arbitrarily
small (i.e., $\mu=\mathcal{U} [-\epsilon,\epsilon]$, $\epsilon>0$).
In Section~\ref{seca-parabolic-version} we show how this sampling
algorithm allows us to extend a classic Monte Carlo scheme of
Muller~\cite{muller1956some},
the so-called ``random walks over spheres'' from the elliptic to the
parabolic setting.

The key idea in this paper is to relate the solution of the obstacle
problem describing the Root barrier with the solution of a nonlinear
integral equation. This general approach dates back to the work of
McKean~\cite{mckean1965appendix}, who showed that the value function
in the pricing problem for a discounted American call option can be
represented in terms of the free boundary function, which itself satisfies
a system of nonlinear integral equations. The question of the uniqueness
of the solution of the integral equation in the context of American
options was resolved by Peskir~\cite{peskir2005american}; see also
the work of Chadam and Chen~\cite{chen2007mathematical}.

Let us finish by mentioning that there have been a number of exciting
recent developments relevant to topics treated in this paper: the
work of Beiglb\"{o}ck and Huesmann~\cite{beiglboeck2013optimal} deriving
the existence of such barriers via optimal transport, the paper of
Galichon, Henry-Labord\`ere, and Touzi who study~(\ref{SEP}) as an
optimal stopping problem~\cite{galichon2011stochastic}, and the work
of Ankirchner, Hobson, and Strack on finite embeddings~\cite
{2013arXiv13063942A,ankirchner2011skorokhod}.

%s2 #&#
\section{The Root barrier as the unique solution of an integral equation}\label{secintegralequatoin}

We begin by recalling classic results on the existence of such barriers.

%de1 #&#
\begin{defn}
\label{defrootbarrier}A closed subset $R$ of $ [0,\infty]\times
[-\infty,\infty]$
is a \emph{Root barrier}~$R$~if:
\begin{longlist}[(2)]
\item[(1)] $ (t,x )\in R$ implies $ (t+h,x )\in R$ $\forall h\geq0$,
\item[(2)] $ (+\infty,x )\in R$ $\forall x\in[-\infty,\infty]$,
\item[(3)] $ (t,\pm\infty)\in R$ $\forall t\in[0,+\infty]$.
\end{longlist}
Given a Root barrier $R$, its \emph{barrier function} $r\dvtx  [-\infty,\infty]\rightarrow[0,\infty]$
is defined as
\[
r (x ):=\inf\bigl\{ t\geq0\dvtx  (t,x )\in R \bigr\},\qquad x\in[-\infty,\infty].
\]
\end{defn}

Note that different barriers can embed the same law. This was resolved
by Loynes by the introduction of regular barriers.

%de2 #&#
\begin{defn}
We say that a barrier $R$, respectively, its barrier function $r$, is
\emph{regular} if $r$ vanishes outside
\[
[x_{-},x_{+} ],
\]
where $x_{+}$ and $x_{-}$ are the first positive, respectively, negative,
zeros\footnote{The first positive zero of some lower-semicontinuous
function $\overline{r}\dvtx  [-\infty,\infty]\rightarrow[0,\infty]$
is at $x_{+}$ if $x_{+}\in[0,\infty]$, $\overline{r} (x_{+} )=0$
and $\overline{r} (x )>0$ for $x\in[0,x_{+} )$.
Similarly for the first negative zero $x_{-}\in[-\infty,0 ]$;
see \cite{MR0292170}, Section~3.} of $r$.
\end{defn}

%th1 #&#
\begin{teo}[(Root \cite{MR0238394}, Loynes \cite{MR0292170}, and Rost \cite{MR0445600})]
\label{teoroot,loynes,rost}Let $\mu$ be a probability measure
on $ (\mathbb{R},\mathcal{B} (\mathbb{R} ) )$ that
has zero mean. Then:
\begin{longlist}[(2)]
\item[(1)] there exists exactly one regular Root barrier $R$ such that $\tau
=\inf\{ t\dvtx  (t,\break B_{t} )\in R \} $
solves (\ref{SEP});
\item[(2)] its barrier function $r (x )=\inf\{ t\dvtx  (t,x )\in R \} $
is a lower semicontinuous function $r\dvtx  [-\infty,\infty]\rightarrow
[0,\infty]$
with $r (\pm\infty)=0$;
\item[(3)] $R= \{ (t,x )\in[0,\infty]\times[-\infty,\infty]\dvtx t\geq r (x )
\} $.
\end{longlist}
Moreover, $\tau$ minimises for every $t\geq0$ the residual expectation
$\mathbb{E} [ (\tilde{\tau}-t )^{+} ]=\int_{t}^{\infty}\mathbb{P}
(\tilde{\tau}>s )\,\dd s$
among all $\tilde{\tau}$ that are solutions of (\ref{SEP}).
\end{teo}

%
%re1 #&#
\begin{rem}
In \cite{MR0238394,MR0292170} the above properties (1)--(3) are only
proved under the additional assumption that $\mu$ is of finite
variance. However, with the help of the PDE representation from \cite
{cox2011root,oberhauserdosreis12} one sees that the finite
variance assumption is unnecessary.
The details may be found in \cite{oberhauserdosreis12}.
%Indeed combining Remark 4.5 in \cite
%{cox2011root} with Lemma 4 in \cite{monroe1972embedding} gives the
%existence of a solution to \eqref{SEP} as a hitting time to a barrier.
%The remaining properties as well as uniqueness can be deduced by
%identification of that barrier with a free boundary as in \cite
%{cox2011root,oberhauserdosreis12}.
\end{rem}

%
%re2 #&#
\begin{rem}
Since the Root barrier $R$ is a closed set, and the process $ (t,B_{t}
)_{t\geq0}$
has continuous trajectories, the representation of $R$ as in point
(3) of Theorem~\ref{teoroot,loynes,rost} above yields
%
%e2.1 #&#
\begin{equation}
\tau\geq r (B_{\tau} ).\label{eqrinequality}
\end{equation}
For example, for $\mu=\frac{1}{2}\delta_{-1}+\frac{1}{2}\delta_{1}$
this is a strict inequality a.s., but in Lemma~\ref{lemrBtau} below
we show that for every atom-free measure, (\ref{eqrinequality}) becomes
an equality. This is intuitive but not completely trivial since it,
for example, also covers the case of singular measures (i.e., not absolutely
continuous with respect to Lebesgue measure but still atom-free) like
Cantor's distribution (devil's staircase).
\end{rem}

%re3 #&#
\begin{rem}
A stopping time $\tau$ minimises the residual expectation if and
only if it minimises for every convex function [wlog $f (0 )=f^{\prime}
(0+ )=0$]
\[
\mathbb{E} \bigl[f (\tau) \bigr]=\int_{0}^{\infty} (
\tau-t )_{+}f^{\prime\prime} (\dd t ).
\]
\end{rem}

Denote the semigroup operator of standard Brownian motion with $
(P_{t}^{B} )$.
The potential kernel is defined as $U^{B}=\int_{0}^{\infty
}P_{t}^{B}\,\dd t$;
that is, $U^{B}$ can be seen as a linear operator on the space of measures
on $ (\mathbb{R},\mathcal{B} (\mathbb{R} ) )$ by
setting $\mu U^{B}=\int_{0}^{\infty}\mu P_{t}^{B}\,\dd t$ which
is of course nothing else than the occupation measure along Brownian
trajectories started with $B_{0}\sim\mu$. If $\mu$ is a signed measure
with $\mu(\mathbb{R} )=0$ and finite first moment, then
the Radon--Nikodym density with respect to the Lebesgue measure is
given as
\[
\frac{\dd \mu\, U^{B}}{\dd x}=-\int\llvert x-y\rrvert\mu
(\dd y ).
\]
Since (in dimension one) Brownian motion is recurrent, $\mu U^{B}$
is infinite if $\mu$ is positive. However, the right-hand side $-\int
_{\mathbb{R}}\llvert x-y\rrvert\mu(\dd y )$
is still well defined for every~$\mu$ that has a finite moment, and
Chacon \cite{MR0501374} demonstrated that this is indeed a very useful
quantity to study hitting times. It plays an essential role for understanding
the dynamics of the Root solution.

%de3 #&#
\begin{defn}[(Potential function)]
Let $\mu$ be a probability measure on $ (\mathbb{R},\mathcal{B} (\mathbb
{R} ) )$
that has a first moment. We define $u_{\mu}\in C (\mathbb{R}, (-\infty,0 ] )$
as
\[
u_{\mu} (x ):=-\int_{\mathbb{R}}\llvert x-y\rrvert\mu(
\dd y )
\]
and call $u_{\mu}$ the \emph{potential function} of the probability
measure $\mu$.
\end{defn}

%s2.1 #&#
\subsection{The barrier function solves an integral equation}

%th2 #&#
\begin{teo}
\label{teointegral}
Denote
\[
g (t,x )=\sqrt{\frac{2t}{\pi}}e^{-x^{2}/(2t)}-\llvert x\rrvert
\operatorname{Erfc} \biggl(\frac{\llvert x\rrvert}{\sqrt{2t}} \biggr
)=\mathbb{E}L_{t}^{x},
\]
where $(L_t^x)_{t,x}$ is the Brownian local time, and let $\mu$ be an
atom-free, zero-mean probability measure on $ (\mathbb{R},\mathcal{B}
(\mathbb{R} ) )$.
Then the regular barrier function $r$ of the Root solution for (\ref{SEP})
solves the nonlinear Volterra integral equation
%
%e2.2 #&#
%e2.3 #&#
\begin{eqnarray}\label{eqintequation}
\qquad u_{\delta_{0}} (x )-u_{\mu} (x )=g \bigl(r (x ),x \bigr)-\int
_{ \{ y\dvtx r (y )<r (x ) \} }g \bigl(r (x )-r (y ),x-y \bigr)\mu(\dd y )
\nonumber\\[-8pt]\\[-12pt]
\eqntext{\forall x\in(-\infty,\infty).}
\end{eqnarray}
\end{teo}

%Let $\mu$ be a probability measure on $\left(\mathbb{R},\mathcal{B}
%\left(\mathbb{R}\right)\right)$
%that has zero mean and no atoms. Then the regular
%barrier function $r$ of the Root solution for (\ref{SEP}) solves
%\[
%u_{\delta_{0}}\left(x\right)-u_{\mu}\left(x\right)=\mathbb{E}\left[L_{r
%\left(x\right)}^{x}\right]-\int_{\left\{ y\dvtx r\left(y\right)<r\left(x
%\right)\right\} }\mathbb{E}\left[L_{r\left(x\right)-r\left(y
%\right)}^{x-y}\right]\mathbb{P}\left(B_{\tau}\in\dd y\right)
%\forall x\in\left(-\infty,\infty\right).
%\]
%Here $\left(L_{t}^{x}\right)_{t,x}$ denotes the local time of Brownian
%motion $\left(B_{t}\right)_{t\geq0}$.
%\end{teo}
%The interest in the above theorem is that only the expected local
%times $\left(L_{t}^{x}\right)_{t,x}$ of the unstopped Brownian motion
%appear, which have an explicit formula. Hence, $r$ solves an explicitly
%given nonlinear integral equation.
%\begin{cor}
%\label{corintequation}
%\end{cor}
We prepare the proof of Theorem~\ref{teointegral} with a lemma:

%le1 #&#
\begin{lem}
\label{lemrBtau} If $\mu$ is atom-free, then $r (B_{\tau} )=\tau$
almost surely.
\end{lem}

\begin{pf}
Since $\mu$ is atom-free, the first positive and negative zeros
cannot be~$0$, that is, $x_{+}>0$ and $x_{-}<0$. We now claim that for all
$ (t,x )$ in the Root barrier $R$,
%
%e2.4 #&#
\begin{equation}
\forall h>0,\forall y\neq x\qquad R\cap[t,t+h )\times(x,y )\neq\varnothing
\label{eqrectangle-1}
\end{equation}
[here $ (y,x )$ should be understood as $ (x,y )$
if $x<y$]. Indeed, assume on the \mbox{contrary} that for some $x$ there
exists $h>0$, $y\neq x$ s.t. $R\cap[r (x ),r (x )+h )\times(x,y )
= \varnothing$.

For simplicity, first assume $0<y<x$ and $r (x )>0$. Then
note that due to lower semicontinuity and Loynes regularity of $r$,
we can find $\underline{r}>0$ and $\delta>0$ such that $r (z )\geq
\underline{r}>0$
for every $z\in(-\delta,x )$. Define $y':=\frac{3y+x}{4}$
and $x':=\frac{3x+y}{4}$, and note that $y<y'<x'<x$. We then have
\begin{eqnarray*}
\mathbb{P} (B_{\tau}=x ) & \geq& \mathbb{P} \Bigl[ \bigl\{
B_{s}\in(-\delta,x ),0\leq s\leq\underline{r} \bigr\} \cap\bigl\{
B_{s}\in(y,x ),\underline{r}\leq s\leq r(x) \bigr\}
\\
& &{} \cap\Bigl\{y<\inf_{r(x)\leq s\leq r(x)+h}B_{s}\leq x \leq\sup
_{r(x)\leq s\leq r(x)+h}B_{s} \Bigr\} \Bigr]
\\
& \geq& \mathbb{P} \bigl[B_{s}\in(-\delta,x ),0\leq s\leq
\underline{r} \bigr]
\\
&&{}\times \inf_{z\in(y',x' )}\mathbb{P}
\bigl[B_{s}\in(y-z,x-z ),0\leq s\leq r(x)-\underline{r} \bigr]
\\
&&{}\times
\inf_{z\in(y',x' )}\mathbb{P} \Bigl[y-z<\inf_{0\leq s\leq h}B_{s}
\leq x-z\leq\sup_{0\leq s\leq h}B_{s} \Bigr]
\\
& > & 0.
\end{eqnarray*}

For the case $r (x )=0$ we have either $x=x_{+},y<x$ or
$x=x_{-},y>x$. In this case an analogous argument works.

Now let $t\mapsto B_{t}\equiv B_{t} (\omega)$ be any continuous
path, and let $t$ be such that $r (B_{t} )=t-\delta<t$.
We claim that this implies that for some $s<t$, $r (B_{s} )<s$.
Indeed, if $B_{t-(\delta/2)}=B_{t}$ we are done; otherwise
by (\ref{eqrectangle-1}), there exists $y$ $\in$ $ (B_{t-(\delta/2)},B_{t} )$
s.t. $r (y )<t-\frac{\delta}{2}$. But then by continuity
of $B$, $B_{s}=y$ for some $s$ $\in$ $ (t-\frac{\delta}{2},t )$,
and this $s$ satisfies $s>r (B_{s} )$.

This argument, together with inequality (\ref{eqrinequality})
and the definition of $\tau$ then imply $r (B_{\tau} )=\tau$.
\end{pf}

Using this, we can now prove Theorem~\ref{teointegral}.

\begin{pf*}{Proof of Theorem~\ref{teointegral}}
Note that by definition of $g$, and since $B_\tau\sim\mu$, the
theorem can be restated as
%
%e2.5 #&#
\begin{eqnarray}\label{eqintwithlocaltime}
u_{\delta_{0}} (x )-u_{\mu} (x )=\mathbb{E} \bigl[L_{r (x )}^{x}
\bigr]-\int_{ \{ y\dvtx r (y )<r (x ) \} }\mathbb{E} \bigl[L_{r (x )-r (y )}^{x-y}
\bigr]\mathbb{P} (B_{\tau}\in\dd y )\nonumber
\\
\eqntext{\forall x\in(-\infty,\infty).}
\end{eqnarray}
Now apply the Tanaka--It\^o formula to the process $ (B_{\tau\wedge t}-x
)_{t\geq0}$
to get
%
%e2.6 #&#
\begin{eqnarray}\label{eqabs}
\mathbb{E} \bigl[\llvert B_{\tau\wedge t}-x\rrvert\bigr] & = & \llvert
x\rrvert
+\mathbb{E} \bigl[L_{t\wedge\tau}^{x} \bigr]\nonumber
\\
& = & \llvert x\rrvert+\mathbb{E} \bigl[L_{t}^{x}+
\bigl(L_{\tau}^{x}-L_{t}^{x}
\bigr)1_{t>\tau} \bigr]
\\
& = & \llvert x\rrvert+g (t,x )-\mathbb{E} \bigl[ \bigl
(L_{t}^{x}-L_{\tau}^{x}
\bigr)1_{t>\tau} \bigr].
\nonumber
\end{eqnarray}
Note that if $\mu$ is atom-free, then $r$ does not have jumps, and
it holds that $\tau=r (B_{\tau} )$ a.s. We use this to transform
the last term into an explicit integral by conditioning%
\footnote{Without loss of generality, we realise Brownian motion on the
canonical Wiener space to justify
the disintegration with the conditional expectation.%
} on $ \{ B_{\tau}\in dy \} $ to see that for all $ (t,x )$
\begin{eqnarray*}
\mathbb{E} \bigl[ \bigl(L_{t}^{x}-L_{\tau}^{x}
\bigr)1_{t>\tau} \bigr] & = & \int_{-\infty}^{\infty}
\mathbb{E} \bigl[ \bigl(L_{t}^{x}-L_{\tau}^{x}
\bigr)1_{t>\tau}\llvert B_{\tau}=y \bigr]\mathbb{P} (B_{\tau}
\in\dd y )
\\
& = & \int_{-\infty}^{\infty}\mathbb{E}
\bigl[ \bigl(L_{t}^{x}-L_{r (y )}^{x}
\bigr)1_{t>r (y )}\llvert B_{\tau}=y \bigr]\mathbb{P} (B_{\tau}
\in\dd y ),
\end{eqnarray*}
where we have used Lemma~\ref{lemrBtau} for the second equality.
If we restrict attention to points $ (r (x ),x )\in R$,
then
\begin{eqnarray*}
\mathbb{E} \bigl[ \bigl(L_{r (x )}^{x}-L_{\tau}^{x}
\bigr)1_{r (x )>\tau} \bigr] & = & \int_{-\infty}^{\infty}
\mathbb{E} \bigl[ \bigl(L_{r (x )}^{x}-L_{r (y )}^{x}
\bigr)1_{r (x )>r (y )}\llvert B_{\tau}=y \bigr]\mathbb{P} (B_{\tau}
\in\dd y )
\\
& = & \int_{ \{ y\dvtx r (y )<r (x ) \} }\mathbb{E} \bigl[L_{r (x
)}^{x}-L_{r (y )}^{x}
\llvert B_{\tau}=y \bigr]\mathbb{P} (B_{\tau}\in\dd y )
\\
& = & \int_{ \{ y\dvtx r (y )<r (x ) \} }\mathbb{E} \bigl[L_{r (x )-r (y )}^{x-y}
\bigr]\mathbb{P} (B_{\tau}\in\dd y ),
\end{eqnarray*}
where for the third equality we have used that Brownian motion is
Markov and that its local time is an additive functional of Brownian
trajectories. Plugging this into (\ref{eqabs}) we see that
\begin{eqnarray*}
\mathbb{E} \bigl[\llvert B_{\tau\wedge r(x)}-x\rrvert\bigr] & = &
\llvert x
\rrvert+\mathbb{E} \bigl[L_{r(x)}^{x} \bigr]-\int
_{ \{ y\dvtx r (y )<r (x ) \} }\mathbb{E} \bigl[L_{r (x )-r (y )}^{x-y} \bigr]
\mathbb{P} (B_{\tau}\in\dd y ).
\end{eqnarray*}
Since $ (r (x ),x )\in R$, the left-hand side multiplied
by $ (-1 )$ equals the potential function of $\mu$, $u_{\mu}$
(see \cite{cox2011root,oberhauserdosreis12} for a proof of this),
and we have derived \eqref{eqintwithlocaltime}.
\end{pf*}

In Section~\ref{subunique} we show that $r$ is not only one
but the unique solution of integral equation (\ref{eqintequation}).
In general it can be hard to numerically solve the integral equation
due to the appearance of the unknown $r$ as an argument in the continuous
integral kernel $g$ as well as the domain of integration. However,
in special cases where more is known about the geometry of $R$, this
can become significantly easier, and in the rest of this article we
focus on measures that lead to symmetric, bounded, and monotone barrier
functions.

%as1 #&#
\begin{assumption}
\label{asnmuhasbarrierfunction}$\mu$ is a zero-mean probability
measure on $ (\mathbb{R},\mathcal{B} (\mathbb{R} ) )$
such that the regular Root barrier solving (\ref{SEP}) is given by
a function $r$ that is symmetric around $0$, continuous, and nonincreasing
on $ [0,\infty]$.
\end{assumption}

The symmetry, boundedness, and especially the monotonicity allows us to
write the integral as an integral with a domain that does not depend
on $r$. This simplifies the numerics needed to solve such integral
equations for the unknown function $r$.

%co1 #&#
\begin{cor}
\label{corintequation-1}Let $\mu$ fulfill Assumption~\ref{asnmuhasbarrierfunction}.
Then $r$ solves the nonlinear Volterra integral equation of the first
kind
%
%e2.7 #&#
%e2.8 #&#
\begin{eqnarray}\label{eqintequation-1}
u_{\delta_{0}} (x )-u_{\mu} (x )
&=& g \bigl(r (x ),x \bigr)\nonumber
\\
&&{} -\int
_{x}^{\infty} \bigl(g \bigl(r (x )-r (y ),x+y \bigr)
+g \bigl(r (x )-r (y ),x-y \bigr) \bigr)\mu(\dd y )
\\
\eqntext{\forall x\in(0, \infty).}
\end{eqnarray}
\end{cor}

\begin{pf}
By assumption on $r$,
\[
\bigl\{ y\dvtx r (y )<r (x ) \bigr\} = (-x,-\infty)\cup(x,\infty),
\]
and by symmetry of the local time in space the statement follows.
\end{pf}

Of course, Assumption~\ref{asnmuhasbarrierfunction} is not too
useful in practice since in general it can be very difficult to deduce
properties
of the geometry of the barrier $R$ from $\mu$. Therefore we provide
in Section~\ref{secbarrierfunction} simple and easy to verify conditions
on $\mu$ that imply Assumption~\ref{asnmuhasbarrierfunction}.

%re4 #&#
\begin{rem}
The solution $\tilde{r}$ of the equation
$u_{\mu} (x )-u_{\delta_{0}} (x )=g (\tilde{r} (x ),x )$
will be a lower bound for the true solution $r$, that is, $\tilde{r} (x
)\leq r (x )$.
Hence a simple inverse problem (or even\vspace*{1pt} a simple ODE after taking
$\frac{\dd}{\dd x}$ if smoothness or $r$ is known) gives a lower bound
for $r$ which often is quite good (e.g., if $\mu=\mathcal{U} [-1,1 ]$).
\end{rem}

%s2.2 #&#
\subsection{The barrier function is the unique solution of the integral equation}
\label{subunique}

We want to find the Root barrier by solving integral
equation (\ref{eqintequation}). Therefore we still need to show
that (\ref{eqintequation}) has a unique solution in a reasonable
class of functions. Unfortunately, there are very few general results
for the uniqueness of such nonlinear integral equations (Volterra's
equation of the first kind); see \cite{linz1985analytical}, Chapter~5.
However, by using the special structure of equation (\ref{eqintequation})
and the connections with viscosity solutions of obstacle PDEs \cite
{oberhauserdosreis12},
we are able to prove uniqueness in the case when $r$ is continuous.
While Theorem~\ref{teointegral} applies to singular distributions
(like the Cantor distribution) we show the uniqueness for solutions of
(\ref{eqintequation})
only for barriers that have a continuous barrier function.

%th3 #&#
\begin{teo}
\label{propintegralrepresentationpde}Let $\mu$ be an atom-free
and zero-mean probability measure on $ (\mathbb{R},\mathcal{B} (\mathbb
{R} ) )$.
If $r\dvtx  (-\infty,\infty)\to[0,\infty]$ is any
continuous function%
\footnote{Note that $r$ is defined as a function taking values that may include
$\infty$; hence $r$ can be continuous, and $r (x )=\infty$
for a $x\in(-\infty,\infty)$ is still possible.%
} that fulfills (\ref{eqintequation}), then
%
%e2.9 #&#
\begin{eqnarray}\label{equrho}
u^{r} (t,x ) &:=& -\int_{-\infty}^{\infty}|y|p (t,x-y
)\,\dd y
\nonumber\\[-8pt]\\[-8pt]
&&{} +\int_{0}^{t}\!\int_{-\infty}^{\infty}1_{\{t\geq r (y )\}}p
(x-y,t-s )\mu(\dd y )\,\dd s\nonumber
\end{eqnarray}
is a continuous viscosity solution with linear growth in the space
variable to
%
%e2.10 #&#
\begin{equation}
\cases{ \min\bigl(u-u_{\mu}, \partial_{t}u-
\frac{1}{2}\partial_{xx}u \bigr) = 0, &\quad on $[0,\infty)
\times\mathbb{R}$,
\vspace*{3pt}\cr
u (t,x ) = -\llvert x\rrvert, &\quad on $\{ 0
\} \times\mathbb{R}$.}\label{eqobstaclepde}
\end{equation}
\end{teo}

\begin{pf}
$u^{r}$ is continuous on $ [0,\infty)\times\mathbb{R}$
and has linear growth in space by standard computations. By defining
\[
Q^{r}:= \bigl\{ (t,x ),t<r (x ) \bigr\}
\]
it is enough to prove:
\begin{longlist}[(3)]
\item[(1)] $\partial_{t}u^{r}-\frac{1}{2}\partial_{xx}u^{r}\geq0$ in viscosity
sense,\vspace*{1.5pt}
\item[(2)] $\partial_{t}u^{r}-\frac{1}{2}\partial_{xx}u^{r}=0$ on $Q^{r}$ in
classical sense,\vspace*{1pt}
\item[(3)] $u^{r} (t,x )\geq u_{\mu} (x )$ on $\mathbb{R}_{+}\times\mathbb{R}$,
and $u^{r} (t,x )=u_{\mu} (x )$ on $ (\mathbb{R}_{+}\times\mathbb{R}
)\setminus Q^{r}$.
\end{longlist}
(1) and (2) are actually true for an arbitrary measurable $r$: indeed,
since $p$ is the fundamental solution to the heat equation, $u^{r}$
solves in a weak (distribution) sense $ (\partial_{t}-\frac
{1}{2}\partial_{xx} )u=1_{\{t\geq r (x )\}}\mu(\dd x )\geq0$,
and the claim follows since distribution (super)solutions to $\partial
_{t}u-\frac{1}{2}\partial_{xx}u=0$
are actually viscosity (super)solutions \cite{ishii1995equivalence}.

It remains to prove point (3). Therefore denote with $p (t,x )=\frac
{1}{\sqrt{2\pi t}}e^{-x^{2}/(2t)}$
the heat kernel. By using Fubini's theorem and that $g (t,x )=\int
_{0}^{t}p (s,x )\,\dd s$,
we immediately see that
\begin{eqnarray*}
& & \int_{0}^{r (x )}\!\int_{-\infty}^{\infty}1_{\{r (x )\geq r (y )\}}p
\bigl(x-y,r (x )-s \bigr)\mu(\dd y )\,\dd s
\\
&&\qquad = \int_{ \{ y\dvtx r (y )<r (x ) \} }g \bigl(r (x )-r (y ),x-y \bigr)\mu(\dd y ).
\end{eqnarray*}
Hence the statement that $r$ solves (\ref{eqintequation}) is equivalent
to the statement
\[
u^{r} \bigl(x,r (x ) \bigr)=u_{\mu} (x ).
\]
Now since $\partial_{xx}u_{\mu}\leq0$, it follows by (2) and comparison
for the heat equation on~$Q^{r}$, that $u^{r}\geq u_{\mu}$ on $Q^{r}$.
To prove that $u^{r}=u_{\mu}$ on $ (\mathbb{R}_{+}\times\mathbb{R}
)\setminus Q^{r}$,
we again use comparison for the heat equation to get that $u^{r}$
is the unique (weak) solution with linear growth to
%
%e2.11 #&#
\begin{equation}
\cases{ \partial_{t}u-\frac{1}{2}
\partial_{xx}u = \mu(\dd x ), &\quad on $(\mathbb{R}_{+}
\times\mathbb{R} )\setminus Q^{r}$,
\vspace*{3pt}\cr
\displaystyle u (t,x ) =
u_{\mu} (x ), &\quad on $ \bigl\{ t=r (x ) \bigr\}$.}\label{eqheatcomplement}
\end{equation}
Note that we use the continuity of $r$ here since it guarantees that
$ (\mathbb{R}_{+}\times\mathbb{R} )\setminus Q^{r}$ is open
and its parabolic boundary is $ \{ t=r (x ) \} $.
\end{pf}

%re5 #&#
\begin{rem}
The representation (\ref{equrho}) is not surprising considering
the classic literature on free boundaries and integral equations cited
in the \hyperref[sec1]{Introduction}. For the Root solution it seems to have been so
far only considered for a special case of the reversed Root (``Rost
barrier'') barrier and derived via pure PDE/nonprobabilistic arguments%
\footnote{We would like to thank Cox for bringing \cite{mcconnell1991two}
to our attention.}~\cite{mcconnell1991two}.
\end{rem}

%
%co2 #&#
\begin{cor}
Let $\mu$ be an atom-free and zero-mean probability measure on $
(\mathbb{R},\mathcal{B} (\mathbb{R} ) )$
that has a continuous barrier function. Then the barrier function
$r$ of the Root solution of (\ref{SEP}) is the unique continuous
function that solves the integral equation (\ref{eqintequation}).
\end{cor}

\begin{pf}
Assume $\overline{r}$ is any other continuous function that solves
(\ref{eqintequation}). Then by Theorem~\ref{propintegralrepresentationpde}
above we know that $u^{\overline{r}}$ and $u^{r}$ both solve the
obstacle PDE (\ref{eqobstaclepde}); hence by the uniqueness result
in \cite{oberhauserdosreis12} they coincide (with $-\mathbb{E} [\llvert
B_{t\wedge\tau}-x\rrvert ]$
where $\tau=\inf\{ t>0\dvtx t\geq r (B_{t} ) \} $).
It follows [by comparing $ (\partial_{t}-\frac{1}{2}\partial_{xx}
)u^{\overline{r}}$
with $ (\partial_{t}-\frac{1}{2}\partial_{xx} )u^{r}$] that
$r (x )=\overline{r} (x )$, $\mu(\dd x )$
a.e., and by continuity and Loynes regularity this implies $r=\overline{r}$.
\end{pf}

%
%re6 #&#
\begin{rem}
The uniqueness result presented here applies to a smaller class of
measures than the class for which integral equation (\ref{eqintequation})
holds. While it covers some cases when the barrier function equals
$\infty$, it does not apply to barriers that arise from singular
measures like the Cantor distribution: while the first two steps of
Theorem~\ref{propintegralrepresentationpde} still hold, we are
not aware of a uniqueness result for the heat PDE~(\ref{eqheatcomplement})
on a complicated (fractal like) domain as $ (\mathbb{R}_{+}\times\mathbb
{R} )\setminus Q^{r}$
(it may no longer be an open set in this case).
\end{rem}

%s3 #&#
\section{Measures with symmetric, continuous and monotone barrier functions}\label{secbarrierfunction}

Assumption~\ref{asnmuhasbarrierfunction}, as introduced in Section~\ref{secintegralequatoin}, is usually not easy to verify for a
given measure $\mu$. It makes a statement about the shape of the
barrier $R$, respectively, $r$, and in general it is very hard to derive
such properties from basic principles. In this section we use the
viscosity methods developed in \cite{oberhauserdosreis12} to show
that simple and easy to verify conditions imply Assumption~\ref
{asnmuhasbarrierfunction}.

%as2 #&#
\begin{assumption}
\label{asnmu}$\mu$ is a symmetric probability measure around $0$
with compact support $ [-k,k ]$ and admits a bounded density
$f$ s.t. $f$ is nondecreasing on~$ [0,k ]$.
\end{assumption}

%re7 #&#
\begin{rem}
If $\mu$ fulfills Assumption~\ref{asnmu}, then $u_{\mu}$ is twice
differentiable on $ (-k,k )$ with
\[
\partial_{xx}u_{\mu}=2f (x ).
\]
\end{rem}

%
%pr1 #&#
\begin{prop}
If $\mu$ fulfills Assumption~\ref{asnmu}, then the corresponding
barrier function $r\dvtx  [-\infty,\infty]\rightarrow[0,\infty]$
is a continuous and nonincreasing function on~$ [0,k ]$.
\end{prop}

\begin{pf}
We first prove the monotonicity. Define $u (t,x )=-\mathbb{E} [\llvert
B_{t\wedge\tau}-x\rrvert ]$.
From \cite{oberhauserdosreis12} it follows that $u$ is the unique
viscosity solution of
%
%e3.1 #&#
\begin{equation}\label{eqpde-1}
\cases{ \min\bigl(u-u_{\mu}, \partial_{t}u-
\frac{1}{2}\partial_{xx}u \bigr) = 0, &\quad on $[0,\infty)
\times[-k,k ]$,
\vspace*{3pt}\cr
\displaystyle u (t,x ) = u_{\delta_{0}} (x ), &\quad on $
\mathbb{R}_+ \times\{ -k,k\} \cup\{ 0 \} \times[-k,k ]$}\hspace*{-25pt}
\end{equation}
and that
%
%e3.2 #&#
\begin{equation}
r (x )=\inf\bigl\{ t\dvtx u (t,x )=u_{\mu} (x ) \bigr\}.\label{eqdef-f}
\end{equation}
We now prove that for any $t\geq0$,
\[
x\mapsto(u-u_{\mu} ) (t,x )
\]
is nonincreasing on $ [0,k ]$ which then implies that $r$
is nonincreasing. Therefore fix a sequence such that $\delta
_{\varepsilon}\to\delta_{0}$
weakly, where $\delta_{\varepsilon}$ has density $\rho^{\varepsilon}$
smooth, symmetric around $0$, decreasing on $\mathbb{R}_{\geq0}$
and support contained in $ [-\varepsilon,\varepsilon]$.
We will consider the functions $u^{\varepsilon}$, unique viscosity
solutions to
%
%e3.3 #&#
\begin{equation}
\cases{ \min\bigl(u^{\varepsilon}-u_{\mu},
\partial_{t}u^{\varepsilon}-\frac{1}{2}\Delta u^{\varepsilon}
\bigr) = 0, &\quad on $[0,\infty)\times[-k,k ]$,
\vspace*{3pt}\cr
\displaystyle u (0,x ) =
u_{\delta_{\varepsilon}} (x ), &\quad on $\displaystyle\mathbb{R}_+
\times\{ -k,k \}
\cup\{ 0 \} \times[-k,k ]$.}\hspace*{-30pt}\label{eqpde-eps}
\end{equation}
Note that since $u_{\delta_{\varepsilon}} (x )\geq u_{\delta_{0}} (x
)-\varepsilon1_{ \{ |x|\leq\varepsilon\}}$,
we have that $u_{\delta_{\varepsilon}}\geq u_{\mu}$ for $\varepsilon$
small enough, and then $u^{\varepsilon}$ admits the representation
%
%e3.4 #&#
\begin{equation}
u^{\varepsilon} (t,x )=-\mathbb{E}_{\delta_{\epsilon}} \bigl[\llvert
B_{t\wedge\tau^{\varepsilon}}-x\rrvert\bigr],
\end{equation}
where $B_{\tau^{\varepsilon}}$ has distribution $\mu$ (for initial
distribution $B_{0}\sim\delta_{\varepsilon}$). The proof now proceeds
in 3 steps.

\begin{longlist}[\textit{Step} 2.]
\item[\textit{Step} 1.] $\partial_{x}u^{\varepsilon}$ exists and is continuous
on $ [0,\infty)\times[-k,k ]$.

For each $t$, $u^{\varepsilon} (t,\cdot)$ has for second
(weak) derivative the measure $\mu_{t\wedge\tau^{\varepsilon}}$,
law of $B_{t\wedge\tau^{\varepsilon}}$. But actually $\mu_{t\wedge\tau
^{\varepsilon}}$
has a bounded density (uniformly in $t\geq0$) since
\begin{eqnarray*}
\mathbb{P}_{\delta_{\varepsilon}} (B_{t\wedge\tau^{\varepsilon}}\in A )
& \leq&
\mathbb{P}_{\delta_{\varepsilon}} (B_{t}\in A )+P (B_{\tau^{\varepsilon
}}\in A )
\\
& \leq& \sup_{t\geq0}\mathbb{P}_{\delta_{\varepsilon}} (B_{t}
\in A )+\mu(A )
\\
& \leq& \bigl(C_{\varepsilon}+\|f\|_{\infty} \bigr)\lambda(A ).
\end{eqnarray*}
Here $\lambda$ is the Lebesgue measure. It follows that $\partial_{x}u$
exists and is continuous in $x$, uniformly in $t$. Joint continuity
then follows easily as in \cite{MR0385326}, Corollary 2.7.

\item[\textit{Step} 2.] $\partial_{x}u^{\varepsilon}\leq\partial_{x}u_{\mu}$
on $ [0,\infty)\times[0,k ]$. Set
\[
D^{+}= \bigl\{ (t,x )\in(0,T ]\times(0,1 )\dvtx  u^{\varepsilon} (t,x
)>u_{\mu} (t,x ) \bigr\}.
\]
We first verify that $w:=\partial_{x}u^{\varepsilon}-\partial_{x}u_{\mu
}\leq0$
on $ ( [0,\infty)\times[0,k ] )\setminus D$:
\begin{itemize}
\item For $x\in[0,k ]$, a direct computation gives $w (0,x )=-2 (\delta
_{\varepsilon}-\mu) [0,x ]$.
Hence
\[
\partial_{x}w (0,x )=-2 \bigl(\rho^{\varepsilon}-f \bigr) (x )
\]
is increasing; that is, $w (0,\cdot)$ is convex, and since
$w (0,k )=w (0,0 )=0$, it follows that $w (0,x )\leq0$, for any $x\in
[0,k ]$.
\item$w (t,0 ) = 0$ since by symmetry $u^{\delta_{\varepsilon}}
(t,x )=u^{\delta_{\varepsilon}} (t,-x )$
(and idem for $u^{\mu}$).
\item On the remaining part $u^{\varepsilon}\equiv u^{\mu}$ so that
$w\equiv0$.
\end{itemize}
Now note that $w$ satisfies
\[
\partial_{t}w-\tfrac{1}{2}\partial_{xx}w=-
\partial_{x}f
\]
(in the distributional sense) on $D^{+}$, and since by assumption
$\partial_{x}f$ is a positive measure, $w$ is a subsolution to the
heat equation on $D^{+}$. Moreover, by step~1 $w$ is continuous
and $w\leq0$ on $\partial D^{+}$, amd hence it follows by the maximum
principle that $w\leq0$ on $D^{+}$ as well.

\item[\textit{Step} 3.] For any $t\geq0$, $x\mapsto(u-u_{\mu} ) (t,x )$
is nonincreasing on $ [0,k ]$.

This is a simple consequence of step~2 and the fact that $u^{\varepsilon
}\to u$
by stability of viscosity solutions.

Hence we get the desired monotonicity of $ (u-u_{\mu} ) (t,\cdot)$
for all $t$, and monotonicity of $r$ follows. It follows that any
discontinuity of $r$ must be of jump-type, but it is obvious that
if $r$ jumps at $x$, then the distribution of $B_{\tau}$ would
have an atom at $x$, which is impossible since $\mu$ has a density.
Hence $r$ is continuous.\quad\qed
\end{longlist}\noqed
\end{pf}

To show that $r (0 )=\sup_{x}r (x )$ is finite
and to provide explicit bounds, we need to make a quantitative
assumption on how fast the mass near $r (0 )$ changes.

%as3 #&#
\begin{assumption}
\label{asnmunear0}$\forall x>0$, $\mu( [-x,x ] )>0$,
and $\exists\eta\in(0,1 )$ s.t.
\[
\sum_{l=0}^{\infty}\eta^{2l} \bigl
\llvert\ln\bigl(\mu\bigl[0,\eta^{l+1}k \bigr] \bigr) \bigr\rrvert<
\infty.
\]
\end{assumption}
%

%re8 #&#
\begin{rem}
A simple family of measures satisfying Assumptions~\ref{asnmu} and~\ref{asnmunear0} is given by
\[
\mu_{k,\alpha} \bigl( [-x,x ] \bigr)= \biggl(\frac{x}{k}
\biggr)^{\alpha},\qquad 0\leq x\leq k,
\]
or any $k>0$, $\alpha\geq1$. In particular, this includes the family
of uniform distributions $\mathcal{U} [-k,k ]$.
\end{rem}

%
%pr2 #&#
\begin{prop}
\label{proprfinite}If $\mu$ fulfills Assumptions~\ref{asnmu}
and~\ref{asnmunear0}, then the corresponding barrier function $r$
is finite on $ [0,k ]$.
\end{prop}

\begin{pf}
Due to the monotonicity and the fact that $\mu$ charges any neighbourhood
of $0$, it is clear that $r (x )$ is finite for any $x>0$.
We now prove $r (0 )<\infty$. First recall that the probability
for Brownian motion to stay in an interval $ (-a,a )$ is
given by
\begin{eqnarray*}
\mathbb{P} \bigl(B_{s}\in(-a,a ),\ \forall 0\leq s\leq T \bigr)&=&
\frac{4}{\pi}\sum_{n=0}^{\infty}
\frac{1}{2n+1}e^{-({(2n+1)^{2}\pi^{2}T})/({8a^{2}})} (-1 ){}^{n}
\\
&\leq&
\frac{4}{\pi}e^{-({\pi^{2}T})/({8a^{2}})};
\end{eqnarray*}
see \cite{feller2008introduction}, Chapter X, Section~5. For any
$0<x<y\leq k$,
we have
\begin{eqnarray*}
\frac{\mu( [-x,x ] )}{\mu( [-y,y ] )} & = & \mathbb{P} \bigl(\llvert
B_{\tau}\rrvert\leq x
|\llvert B_{\tau}\rrvert\leq y \bigr)
\\
& \leq& \mathbb{P} \Bigl(\sup_{r(y)\leq s\leq r(x)}\llvert
B_{s}-B_{r(y)}\rrvert\leq2y \Bigr)
\\
& \leq& \frac{4}{\pi}e^{-({\pi^{2} (r(x)-r(y) )})/({32y^{2}})}.
\end{eqnarray*}
This can be rewritten as
%
%e3.5 #&#
\begin{equation}
r (x )\leq r (y )+\frac{32y^{2}}{\pi^{2}} \biggl(\ln\biggl(\frac{4}{\pi}
\biggr)+ \biggl\llvert\ln\biggl(\frac{\mu( [0,x ] )}{\mu( [0,y ] )}
\biggr) \biggr\rrvert
\biggr).\label{eqbnd-r}
\end{equation}
Now fix $0<\eta<1$ from Assumption~\ref{asnmunear0}. Taking
successively $ (x,y )= (\eta^{l+1}k,\eta^{l}k )$
in (\ref{eqbnd-r}) and summing, we get
\[
r \bigl(\eta^{r+1}k \bigr)\leq\frac{32k^{2}}{\pi^{2}}\sum
_{l=0}^{r}\eta^{2l} \biggl(\ln\biggl(
\frac{4}{\pi} \biggr)+ \biggl\llvert\ln\biggl(\frac{\mu( [0,\eta
^{l+1}k ] )}{\mu( [0,\eta^{l}k ] )} \biggr)
\biggr\rrvert\biggr).
\]
It only remains to let $l\to\infty$, and we finally obtain $r (0^{+}
)<\infty$.
Putting the above together gives us the desired implication.
\end{pf}

%
%co3 #&#
\begin{cor}
If $\mu$ fulfills Assumptions~\ref{asnmu}~and~\ref{asnmunear0}
then $\mu$ fulfills Assumption~\ref{asnmuhasbarrierfunction}.
\end{cor}

The above proofs show much more about $r$ in the sense that they can
give an explicit upper and lower bound on $\sup_{x\in\mathbb{R}}r (x
)=r (0 )$.
For example, for the special case of $\mu=\mathcal{U} [-1,1 ]$
that we are interested in for our Monte Carlo application one easily
derives the following statement.

%co4 #&#
\begin{cor}
\label{corboundsuniform}Let $\mu$ be the uniform distribution
on $ [-1,1 ]$. Then
\[
r(0)\in\biggl[\frac{\pi}{8},\frac{32}{\pi^{2}}\inf_{\eta\in(0,1 )}
\frac{\ln({4}/(\pi\eta) )}{1-\eta^{2}} \biggr].
\]
\end{cor}

%
%f1 #&#
\begin{figure}[t]

\includegraphics{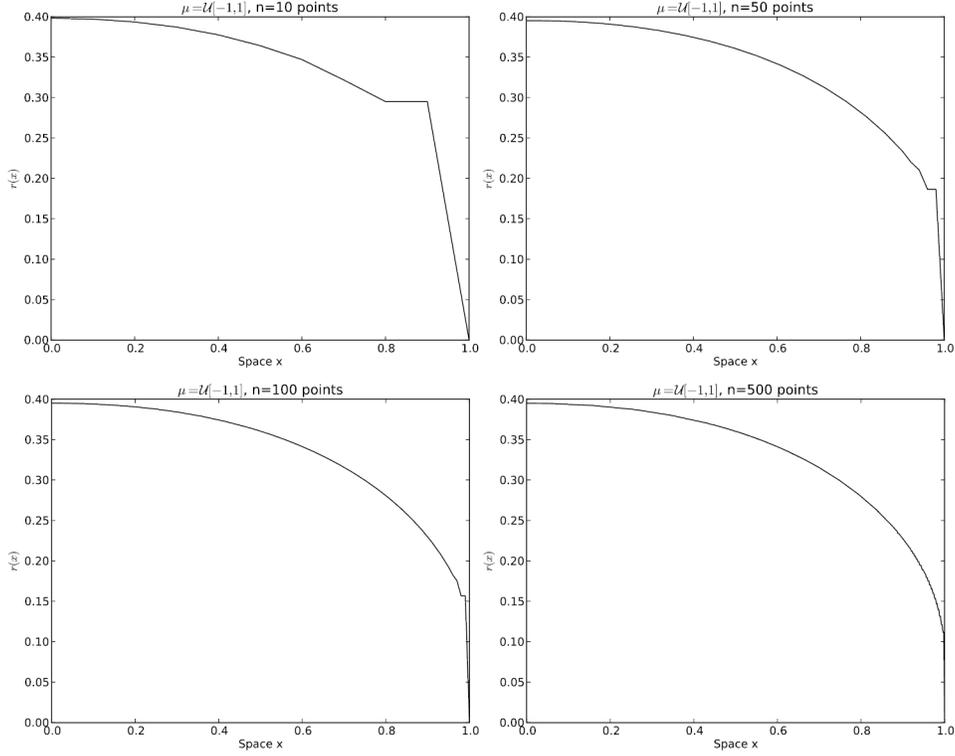}

\caption{The Root barrier for $\mu=\mathcal{U} [-1,1 ]$.
The above figures were produced with the forward Euler scheme implemented
in Python (SciPy \cite{python,Hunter2007}). The integral equation
is stable in the sense that already with only $10$ discretization
points the approximation is fairly accurate away from $x=1$. With
$n=500$ points the program finishes in less than $3$ seconds on
a standard laptop (Intel i5-3210M, 3.10~GHz, 3~MB L3, 1600~MHz FSB,
8~GB DDR3 RAM).}\label{figintegraluniform}
\end{figure}

\begin{pf}
Since $\tau\leq\sup_{x}r (x )=r (0 )$ we have
\[
\mathbb{E} \bigl[\llvert B_{\tau}\rrvert\bigr]\leq\mathbb{E} \bigl[
\llvert B_{r (0 )}\rrvert\bigr].
\]
Using $B_{\tau}\sim\mathcal{U} [-1,1 ]$ and a simple calculation
this becomes
\[
\frac{1}{2}\leq\sqrt{\frac{2}{\pi}r (0 )}
\]
which immediately gives the lower bound. The upper bound follows
from the proof of Proposition~\ref{proprfinite} since in this
case
\[
\frac{\mu( [0,\eta^{r+1}k ] )}{\mu( [0,\eta^{r}k ] )}=\eta.
\]\upqed
\end{pf}

%
%re9 #&#
\begin{rem}
Numerics given in the next section and Figure~\ref{figintegraluniform}
show that this lower bound is actually very good ($\frac{\pi
}{8}=0.392\ldots$)
but that the upper bound $\frac{32}{\pi^{2}}\inf_{\eta\in(0,1 )}\frac
{\ln(4/(\pi\eta) )}{1-\eta^{2}}=3.774\ldots$
is not.

%re10 #&#
\begin{rem}
It is interesting to compare our Proposition~\ref{proprfinite} to the
results of Ankirchner and Strack \cite{ankirchner2011skorokhod}. On one
hand, they obtain a general necessary condition for a bounded time
embedding to exist, namely
\[
\sup_{x \in\operatorname{supp}(\mu)} \limsup_{\varepsilon\downarrow0}
\varepsilon^2 \bigl\llvert\ln\bigl( \mu[x-\varepsilon, x +
\varepsilon] \bigr) \bigr\rrvert< \infty,
\]
where we recognise the term in the series from Assumption~\ref{asnmunear0}.
On the other hand, they also study an embedding
due to Bass and obtain sufficient conditions under which the associated
stopping time $\tau^{B}$ is bounded. Note that any almost sure bound
on $\tau^{B}$ implies the same bound for the Root stopping time $\tau^{R}$
(and hence the barrier function $r$), since $\tau^{R}$ minimises
$\mathbb{E} (\tau-t ){}_{+}$ for all $t\geq0$. In fact, one can check
that under Assumption~\ref{asnmu}, the sufficient conditions given in
\cite{ankirchner2011skorokhod} all imply our Assumption~\ref
{asnmunear0} (of course, this does not mean that their results are a
corollary of ours, since they deal with general measures while we only
have to check the behaviour around the point $0$). In addition, the
upper bounds obtained in \cite{ankirchner2011skorokhod} are sometimes
sharper. For instance,
we could deduce from their results the upper bound $r (0 )\leq\frac
{2}{\pi}=0.636\ldots$
for $\mu=\mathcal{U} [-1,1 ]$; that is, without running numerics
we already know that $\sup_{x}r (x )\in[\frac{\pi}{8},\frac{2}{\pi} ]$.
\end{rem}
\end{rem}

%s3.1 #&#
\subsection{Numerics for the integral equation}

Due to the importance of such an integral equation in engineering and
physics, there is an abundance of literature treating numerics; see
\cite{linz1985analytical} and the reference therein. We therefore
do not discuss proofs of convergence, etc. Instead we give a simple example
that demonstrates that already the arguably simplest scheme, a forward
Euler discretisation,
provides a very fast way to solve the integral equation.

To calculate $r$ for a given $\mu$ with $\operatorname{supp}= [-k,k ]$
and density $f$, fix $n\in\mathbb{N}$, set $h=\frac{k}{n}$, and
for every $i\in\{ 1,\ldots, n \} $ denote with $r_{i}$
the approximation to $r (ih )$. Then we know that $r_{n}=0$,
and (starting with $i=n-1$) we can solve recursively the discretised
nonlinear equation for $r_{i}$,
\begin{eqnarray*}
&& u_{\mu} (ih )-u_{\delta} (ih )
\\
&&\qquad =g (r_{i},ih )-\sum
_{j=i+1}^{n} \bigl(g \bigl(r_{i}-r_{j},
(i-j )h \bigr)+g \bigl(r_{i}-r_{j}, (i+j )h \bigr) \bigr)f
(jh ).
\end{eqnarray*}
%

%s4 #&#
\section{Generating bounded Brownian time--space increments}\label{secgenerating-bounded-brownian}

As an application of the previous sections we now return to the approach
pointed out in the \hyperref[sec1]{Introduction}: that an intelligent choice of $\mu$
can lead to an efficient procedure to sample from Brownian trajectories.

%co5 #&#
\begin{cor}
\label{corRootBM}There exists a continuous bounded function
\[
r\in C_{b} \bigl( [-1,1 ],\mathbb{R} \bigr)\qquad\mbox{with }r (x )=r
(-x )\geq0\mbox{ and } r (1 )=r (-1 )=0
\]

which is decreasing on $ [0,1 ]$ such that:
\begin{longlist}[(2)]
\item[(1)] if $B$ is Brownian motion carried on a probability space $ (\Omega,\mathcal{F},\mathcal{F}_{t},\mathbb{P} )$
satisfying the usual conditions,
\item[(2)] and the \textup{\emph{sequence of}} stopping times $\tau= (\tau
_{k} )_{k\geq0}$
is defined as
\[
\tau_{0}=0\quad\mbox{and}\quad\tau_{k+1}=
\tau_{k}+\inf\bigl\{ \Delta\dvtx \Delta\geq r (B_{\tau_{k}+\Delta}-B_{\tau_{k}}
) \bigr\}
\]
[i.e., $\tau_{1}$ is the exit time from $R= \{ (t,x )\dvtx t\leq r (x ) \} $],
\end{longlist}
then the following properties hold:
\begin{longlist}[(2)]
\item[(1)] if\vspace*{1.5pt} $ (U_{k} )_{k\geq1}$ is a sequence of i.i.d. random
variables carried on a probability space $ (\Omega^{\mathrm{sim}},\mathcal
{F}^{\mathrm{sim}},\mathbb{P}^{\mathrm{sim}} )$,
each uniformly distributed on $ [-1,1 ]$, $U_{1}\sim\mathcal{U} [-1,1 ]$,
then
\[
(\tau_{k+1}-\tau_{k},B_{\tau_{k+1}}-B_{\tau_{k}}
)_{k\geq0}\stackrel{\mathrm{Law}} {=} \bigl(r (U_{k}
),U_{k} \bigr)_{k\geq0},
\]

\item[(2)] $\llvert\tau_{k+1}^{\epsilon}-\tau_{k}^{\epsilon}\rrvert\leq r
(0 )<\infty$
and \textup{$\sup_{t\in[\tau_{k},\tau_{k+1} ]}\llvert B_{t}-B_{\tau
_{k}}\rrvert\leq2$}
for every $k\geq0$.
\end{longlist}
Moreover, the function $r$ is the unique continuous solution of the
integral equation
%
%e4.1 #&#
\begin{eqnarray}
\frac{x^{2}+1}{2}-x&=&g \bigl(r (x ),x \bigr)\nonumber
\\
&&{} -\frac{1}{2}\int
_{x}^{1} \bigl(g \bigl(r (x )-r (y ),x-y \bigr)+g
\bigl(r (x )-r (y ),x+y \bigr) \bigr)\,\dd y\nonumber
\\
\eqntext{\forall x\in[0,1 ],}
\end{eqnarray}
where
\[
g (t,x )=\mathbb{E}L_{t}^{x}=\sqrt{\frac{2t}{\pi}}e^{-x^{2}/(2t)}-
\llvert x\rrvert\operatorname{Erfc} \biggl(\frac{\llvert x\rrvert
}{\sqrt{2t}} \biggr).
\]
\end{cor}

\begin{pf}
This follows directly from Theorem~\ref{teointegral} and Markovianity
of Brownian motion.
\end{pf}

\begin{figure}[b]

\includegraphics{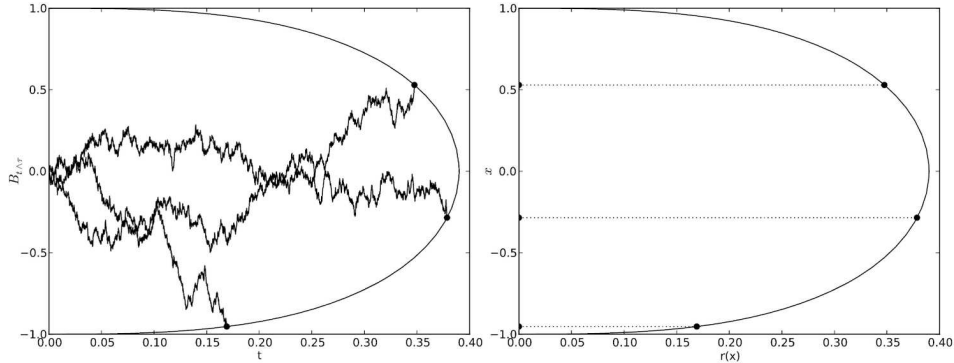}

\caption{The plot on the left shows three Brownian trajectories that
were stopped after hitting the Root barrier for $\mu=\mathcal{U} [-1,1
]$. The plot on the right is the same but with the trajectories removed
and the hitting points of the Root barrier projected back to $\mathbb{R}$.}\label{fig2}
\end{figure}

\begin{figure}[t]

\includegraphics{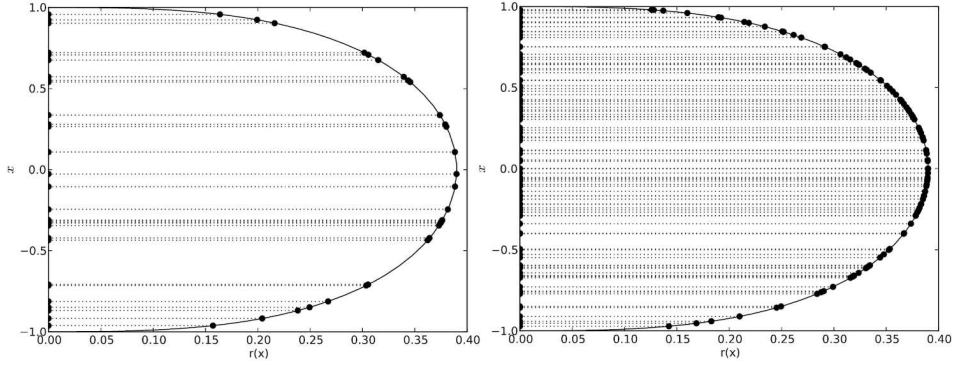}

\caption{Similarly to the above, both plots were drawn by using 30,
respectively, 100, samples from a Brownian motion.
We see can start to see that the projected points follow a $\mathcal{U}
[-1,1 ]$ distribution.}\label{fig3}
\end{figure}

We refer to Figures \ref{fig2} and \ref{fig3} below for some examples of uniformly distributed space increments obtained by the above procedure.

%
%f2 #&#

%
%f3 #&#

Using Brownian scaling one immediately gets:

%co6 #&#
\begin{cor}
\label{corRootBMscaled}If we fix $\epsilon>0$ and replace in the
above the sequence $\tau= (\tau_{k} )$ by $\tau^{\epsilon}= (\tau
_{k}^{\epsilon} )$
defined as
\[
\tau_{0}^{\epsilon}=0\quad\mbox{and}\quad\tau_{k+1}^{\epsilon}=
\tau_{k}^{\epsilon}+\inf\biggl\{ \Delta\dvtx \Delta\geq
\epsilon^{2}r \biggl(\frac{B_{\tau_{k}^{\epsilon}+\Delta}-B_{\tau
_{k}^{\epsilon}}}{\epsilon} \biggr) \biggr\},
\]

then the following properties hold:
\begin{longlist}[(2)]
\item[(1)] for a sequence $ (U_{k} )_{k\geq1}$ of i.i.d. random variables
carried on a probabilty space $ (\Omega^{\mathrm{sim}},\mathcal{F}^{\mathrm{sim}},\mathbb
{P}^{\mathrm{sim}} )$,
each uniformly distributed on $ [-1,1 ]$, $U_{1}\sim\mathcal{U} [-1,1 ]$ we
have
\[
\bigl(\tau_{k+1}^{\epsilon}-\tau_{k}^{\epsilon},B_{\tau_{k+1}^{\epsilon
}}-B_{\tau_{k}^{\epsilon}}
\bigr)_{k\geq1}\stackrel{\mathrm{Law}} {=} \bigl(\epsilon^{2}r
(U_{k} ),\epsilon U_{k} \bigr)_{k\geq1},
\]

\item[(2)] $\llvert\tau_{k+1}^{\epsilon}-\tau_{k}^{\epsilon}\rrvert
<\epsilon^{2}r (0 )$
and $\sup_{t\in[\tau_{k}^{\epsilon},\tau_{k+1}^{\epsilon} ]}\llvert
B_{t}-B_{\tau_{k}^{\epsilon}}\rrvert\leq2\epsilon$
for every $k\geq0$.
\end{longlist}
\end{cor}

The interest in above statement is to simulate time--space Brownian
motion $t\mapsto(t,B_{t} )$ on a computer in an easy and
efficient way: to sample one increment we only need to generate one
uniformly distributed random variable $U$ and evaluate the function
$r$ at $U$ to match in law the increment of the time--space process
$ (\tau_{k+1}-\tau_{k},B_{\tau_{k+1}}-B_{\tau_{k}} )$. In pseudo code it reads Algorithm~\ref{alg1}.

\begin{algorithm}%[b]
\caption{Generate a Brownian increment from $\mathcal{U} [-1,1 ]$}\label{alg1}
\begin{algorithmic}[1]
\Function{SampleBMincrement}{$\epsilon$}
\State$U \gets\mathcal{U} [-1,1 ]$
\State$\Delta B \gets\epsilon* U$
\State$\Delta t \gets\epsilon^2*r(U)$
\State\Return$(\Delta t,\Delta B)$
\EndFunction
\end{algorithmic}
\end{algorithm}

Contrast this with standard methods where the time step is deterministic,
but a normally distributed space increment is simulated by transformations
of (several) uniformly distributed random variables and table look-ups
(e.g., via the Box--Muller transform, the Ziggurat algorithm,
the Marsaglia polar method, etc.).

On the other hand, the function $r$ in above statement is not given
by an explicit analytic expression. However, the integral equation
can be solved with great precision, and this computation needs to be
done only once, then stored in a table (possibly after spline interpolation,
etc.), that is, evaluating $r$ at a point amounts to a table look-up.

The most attractive feature of the above algorithm is that one can \emph{fix
at every step a deterministic bound} \emph{on the space and time increments,
and both resulting increments are trivial to simulate}. In the next
section we demonstrate this advantage on a short and simple but nontrivial
example: a Monte Carlo simulation with adaptive step size applied
to parabolic PDEs. The deterministic control over time--space increments
allows us to make very big steps without leaving the time--space domain
which leads to a very fast algorithm.

%s5 #&#
\section{A parabolic version of Muller's random walk over spheres}\label{seca-parabolic-version}

The use of exit times from a domain to simulate Brownian motion is
classic and goes back to Muller in 1956 who used the uniform exit
distribution of Brownian motion from a sphere to calculate elliptic
PDEs (the so-called ``random walk on touching spheres'') of the form
\[
\cases{ \frac{1}{2}\Delta u = 0, &\quad on $\mathcal{D}$,
\vspace*{3pt}\cr
u (x ) = g
(x ), &\quad on $\partial\mathcal{D}$,}
\]
where $\mathcal{D}$ is a domain in $\mathbb{R}^{n}$ via the Monte Carlo
approximations to $u (x )=\mathbb{E}_{t,x} [g (B_{\tau^{\mathcal{D}}} ) ]$.
Here $\tau^{\mathcal{D}}$ denotes the exit time of $B$ from $\mathcal{D}$.
The attraction of this approach is that in every step one can choose
the diameter of the sphere arbitrarily big, subject only to not intersecting
$\partial\mathcal{D}$ before one samples the Brownian increment.
These give big Brownian increments that lead to a very fast algorithm.
To make this work for a parabolic PDE
%
%e5.1 #&#
\begin{equation}
\cases{ \partial_{t}u+\frac{1}{2}\Delta u = 0, &\quad on $
\mathcal{D}$,
\vspace*{3pt}\cr
u (t,x ) = g (t,x ), &\quad on $\mathcal{PD}$}\label{eqpPDE}
\end{equation}
(here we denote $\mathcal{D}=\bigcup_{t\geq0} \{ t \} \times
D_{t}\subset[0,\infty)\times\mathbb{R}$
and the parabolic boundary $\mathcal{P}\mathcal{D}=\partial\mathcal
{D}\setminus( \{ 0 \} \times D_{0} )$),
it is necessary to additionally sample the distribution of the exit
time from the sphere. While analytic expressions are known, it is
not efficient to simulate. This has been pointed out by many authors
and the work of Milstein and Tretyakov, Deaconu and Hermann, Deaconu,
Lejay, and Zein \mbox{\cite{MR1722281,deaconeherrman2013,MR2673769}}, proposes
the use of exit times of time--space Brownian motion from other shapes
than spheres. The approach which is closest to the one presented here
is the random ``walk over moving spheres'' (WoMS) introduced in \cite
{deaconeherrman2013}.
In the short section below we show that the Root solution gives another
way to construct such a random walk. It is optimal among all such
approaches \cite{MR1722281,deaconeherrman2013,MR2673769} in the sense
that one samples simply from the uniform distribution. A (theoretical)
disadvantage is that the barrier $r$ is not known in explicit form
and has to be stored as a table look-up, though the results from the
previous sections show that this can be done quite easily.

%s5.1 #&#
\subsection{A random walk over Root barriers}

We introduce here a Monte Carlo scheme to calculate the solution
of the parabolic PDE (\ref{eqpPDE}). To avoid technicalities we
assume the boundary is smooth.

%as4 #&#
\begin{assumption}
\label{assumption}The space--time domain is of the form
\[
\mathcal{D}=\bigcup_{t\in(0,T )} \{ t \} \times
(a_{t},b_{t} ),
\]
where $T\in(0,\infty)$ is fixed, $a,b\in C^{1} ( (0,T ),\mathbb{R} )$,
and $a_{t}<b_{t}$ on $ (0,T )$. In addition~$g$ is assumed
to be regular enough so that the solution $u$ to (\ref{eqpPDE})
satisfies
\[
\bigl|u (t,x )-u (s,y )\bigr|\leq\llvert u\rrvert_{\mathrm{Lip}} \bigl(|t-s|^{1/2}+|x-y|
\bigr)\qquad \forall(t,x ), (s,y )\in\mathcal{D}
\]
for some constant $\llvert u\rrvert_{\mathrm{Lip}}<\infty$; see, for example,
\cite{lieberman1996second}
for several standard conditions guaranteeing this.
\end{assumption}
%

%de4 #&#
\begin{defn}
The parabolic distance to the boundary $\mathcal{D}$ is defined as
\[
d_{\mathcal{D}} (t,x )=\min(x-a_{t},x-b_{t},\sqrt{T-t}
).
\]
For $\delta>0$ define $\mathcal{D}_{\delta}$ as
\[
\mathcal{D}_{\delta}= \bigl\{ (t,x )\in\mathcal{D}\dvtx d (t,x )\leq\delta
\bigr\}.
\]
\end{defn}

%
%re11 #&#
\begin{rem}
Since $a,b$ are Lipschitz, one can find a function $\rho=\rho(t,x )$
such that:
\begin{itemize}
\item$c.d_{\mathcal{D}} (t,x )\leq\rho(t,x )\leq d_{\mathcal{D}} (t,x )$
for some constant $c>0$,
\item$\forall(t,x )\in\mathcal{D}$ we have $B_{t,x}^{\rho(t,x
)}\subset\overline{\mathcal{D}}$.
\end{itemize}
\end{rem}

%
%de5 #&#
\begin{defn}
Denote $r$ the barrier function associated with $\mu=\mathcal{U} [-1,1 ]$
and with $R_{t,x}^{\epsilon}$ its Root barrier around $ (t,x )$
after scaling with some $\epsilon>0$, that is,
\[
R_{t,x}^{\epsilon}= \bigl\{ \bigl(t+\epsilon^{2}s,x+
\epsilon y \bigr)\dvtx s\geq r (y ) \bigr\}.
\]
\end{defn}

We now introduce a Markov chain that is easy to generate on a computer.
The motivation is the following: fix a point $ (t,x )\in\mathcal
{D}\setminus\mathcal{D}_{\delta}$,
and consider the Root barrier $R_{t,x}^{\rho(t,x )}$. From
the very definition of $\rho(t,x )$, it follows that a Brownian
motion started at $ (t,x )$ will not have left the domain
$\mathcal{D}$ before it leaves $R_{t,x}^{\rho(t,x )}$.
We now record the exit time and position of $B$ from $R_{t,x}^{\rho
(t,x )}$,
and Corollary~\ref{corRootBMscaled} tells us that the distribution
of this time--space increment is $ (\rho^{2} (t,x )r (U ),\rho(t,x )U )$
for $U\sim\mathcal{U} [-1,1 ]$. If this first step puts
us into $\mathcal{D}_{\delta}$, we stop. Otherwise we carry out the
same procedure again, but now starting at $ (t+\rho^{2} (t,x )r (U
),x+\rho(t,x )U )$.

%
%de6 #&#
\begin{defn}
For every $ (t,x )\in\mathcal{D}$ define a Markov chain
\[
M^{t,x,\delta}= \bigl(\tau_{k}^{t,x,\delta},M_{k}^{t,x,\delta}
\bigr)_{k\geq1}= (\tau_{k},M_{k} )_{k\geq1}
\]
and a stopping time $\nu=\nu^{t,x,\delta}$ (if the context is clear,
we do not write the superscripts $t,x,\delta$) recursively as follows:
\[
(\tau_{0},M_{0} )  =  (t,x )
\]
and
\begin{eqnarray*}
&& (\tau_{k+1},M_{k+1} )
\\
&&\qquad  =  \cases{ \displaystyle\bigl(
\rho^{2} (\tau_{k},M_{k} )r (U_{k}
),M_{k}+\rho(\tau_{k},M_{k} )U_{k}
\bigr), &\quad if $\displaystyle(\tau_{k},M_{k} )\in
\mathcal{D}/\mathcal{D}_{\delta}$,
\vspace*{3pt}\cr
\displaystyle(\tau_{k},M_{\tau_{k}}
), &\quad if $\displaystyle(\tau_{k},M_{k} )\in
\mathcal{D}_{\delta}$.}
\end{eqnarray*}
Further denote $\nu=\inf\{ k\dvtx  (\tau_{k},M_{k} )\in\mathcal{D}_{\delta}
\} $
and
\[
\bigl(\nu^{\mathcal{D}},M_{\nu}^{\mathcal{D}} \bigr)= \cases{
\displaystyle(\nu,a_{\nu} ), &\quad if $\displaystyle d_{\mathcal{D}}
(\nu,M_{\nu} )=a_{\nu}-M_{\nu}$,
\vspace*{3pt}\cr
\displaystyle
( \nu,b_{\nu} ), &\quad if $\displaystyle d_{\mathcal{D}} (
\nu,M_{\nu} )=M_{\nu}-b_{\nu}$,
\vspace*{3pt}\cr
\displaystyle
(T,M_{\nu} ), &\quad otherwise.}
\]

\begin{algorithm}%[b]
\caption{Random walk over Root barriers}\label{alg2}
\begin{algorithmic}[1]
\Function{RootMonteCarlo}{$t,x,samples}$
\State$u\gets0$
\For{$i\gets1,samples$}
\State$(\tau,B) \gets(t,x)$
\While{$\rho(\tau,B)>\delta$}
\State$(\Delta\tau,\Delta B) \gets SampleBMincrement(\rho(\tau,B))$
\State$(\tau,B) \gets(\tau+\Delta\tau,B+\Delta B)$
\EndWhile
\State$u\gets u+g(\tau,B)$
\EndFor
\State$u\gets u/samples$
\State\Return$u$
\EndFunction
\end{algorithmic}
\end{algorithm}

Put simply, once our Markov chain enters $\mathcal{D}_{\delta}$, we
stop it, and $ (\nu^{\mathcal{D}},M_{\nu}^{\mathcal{D}} )$
then records the nearest point on the boundary. This very easy to
implement and spelled out in pseudocode it reads as Algorithm~\ref{alg2}.
\end{defn}

By construction of the Markov chain, it is clear that each sample
trajectory does not contribute an error bigger than $\delta$. The
more interesting question is how many steps the chain makes on average
before leaving $\mathcal{D}_{\delta}$. As in Muller's elliptic version~\cite{muller1956some}, the average number of steps only grows proportionally
to~$\log\frac{1}{\delta}$.

%th4 #&#
\begin{teo}
If Assumption~\ref{assumption} holds, then there exists a unique solution
$u$ in the class \textup{$C^{1,2} (\mathcal{D},\mathbb{R} )\cap C
(\overline{\mathcal{D}},\mathbb{R} )$}
that solves (\ref{eqpPDE}). Moreover, there exist constants $c_{1}$,
$c_{2}$, $\delta_{0}$
such that for every $\delta\in(0,\delta_{0} )$ one has
\[
\bigl\llvert\mathbb{E}_{t,x} \bigl[g \bigl(\tau_{\nu},M_{\nu}^{\mathcal{D}}
\bigr) \bigr]-u (t,x ) \bigr\rrvert\leq c_{1}\delta.
\]
The number of steps $\nu$ is finite a.s., and for all $ (t,x )\in
\mathcal{D}\setminus\mathcal{D}_{\delta}$,
\[
\mathbb{E}_{t,x} [\nu]\leq c_{2} \bigl(1+\log(1/\delta)
\bigr).
\]
\end{teo}

\begin{pf}
Under the above assumptions on $g$ and $\mathcal{D}$, the existence
of a unique classical solution to (\ref{eqpPDE}) and the Feynman--Kac
representation
%
%e5.2 #&#
\begin{eqnarray}
u (t,x )=\mathbb{E} \bigl[g \bigl(\sigma^{t,x}\wedge
T,B_{\sigma^{t,x}\wedge T}^{t,x} \bigr) \bigr]\nonumber
\\
\eqntext{\mbox{where }
\sigma^{t,x}=\inf\bigl\{ s>t\dvtx B_{s}^{t,x}\notin
(a_{s},b_{s} ) \bigr\},}
\end{eqnarray}
and $B^{t,x}$ denotes a Brownian motion started at $x$ at time $t$
follows from the standard results; see, for example, \cite
{lieberman1996second}, Theorems 5.9, 5.10, 6.45,
and for the Feynman--Kac verification, \cite{costantini2006boundary},
Appendix~B.
Write
\begin{eqnarray*}
&& \mathbb{E}_{t,x} \bigl[g \bigl(\tau_{\nu},M_{\nu}^{\mathcal{D}}
\bigr) \bigr]-u (t,x )
\\
&&\qquad  = \mathbb{E}_{t,x} \bigl[u \bigl(
\tau_{\nu},M_{\nu}^{\mathcal{D}} \bigr) \bigr]-
\mathbb{E}_{t,x} \bigl[u (\tau_{\nu},M_{\nu} ) \bigr]
 +\mathbb{E}_{t,x} \bigl[u (\tau_{\nu},M_{\nu} )
\bigr]-u (t,x ),
\end{eqnarray*}
and note that the first difference on the right-hand side is bounded
by $\llvert u\rrvert_{\mathrm{Lip}}\delta$. The second difference on the
right-hand side vanishes since by construction of the Markov chain, we
have $ (\tau_{\nu},M_{\nu} )\stackrel{\mathrm{Law}}{=} (\tau_{\nu
},B_{\tau_{\nu}} )$,
and $u$ is space--time harmonic on~$\mathcal{D}$. To estimate the
number of steps, we start with an idea similar to that in \mbox{\cite{MR1489185,muller1956some}}
but then argue via PDE comparison. This allows us to give a short
proof. For $v$ a bounded measurable function on \emph{$\mathcal{D}$},
define
\begin{eqnarray*}
Pv (t,x ) & =&\mathbb{E}_{t,x} \bigl[v \bigl(\tau^{t,x},B_{\tau^{t,x}}
\bigr) \bigr],
\end{eqnarray*}
where $\tau^{t,x}$ is the first exit time from $R_{t,x}^{\rho(t,x )}$.
We denote the expected number of steps with $n (t,x )=\mathbb{E}_{t,x}
[\nu]$.
It is then the unique solution to the equation
%
%e5.3 #&#
\begin{equation}
\cases{ n-Pn=1, &\quad in $\mathcal{D}\setminus\mathcal{D}_{\delta}$,
\vspace*{3pt}\cr
n=0, &\quad in $\mathcal{D}_{\delta}$.}\label{eqPn}
\end{equation}
To obtain an upper bound on $n$ it is enough to obtain supersolutions
to the above equation. Note that if $v$ is $\mathcal{C}^{1,2} (\overline
{\mathcal{D}} )$,
by It\^{o}'s formula we actually have
%
%e5.4 #&#
\begin{eqnarray}
Pv (t,x ) & =& v (t,x )+\mathbb{E}_{t,x} \biggl[\int_{t}^{\tau^{t,x}}
\biggl(\partial_{t}+\frac{1}{2}\partial_{xx} \biggr)v
(s,B_{s} )\,\dd s \biggr].\label{eqPVIto}
\end{eqnarray}
Now take
\begin{eqnarray*}
v^{1} (t,x ) & = & \log(x-a_{t}+\delta)+\log
(b_{t}-x+\delta)+\tfrac{1}{2}\log\bigl(T-t+\delta^{2}
\bigr),
\end{eqnarray*}
and direct computation shows that for small enough $\eta>0$ (not
depending on $\delta$, assuming if necessary $\delta$ smaller than
some suitable $\delta_{0}$),
\begin{eqnarray*}
&& \biggl(\partial_{t}+\frac{1}{2}\partial_{xx}
\biggr)v^{1} (t,x )
\\
&&\qquad  =  -\frac{1}{2} \biggl(\frac{1}{|x-a_{t}+\delta|^{2}}+
\frac{1}{|b_{t}-x+\delta|^{2}}+\frac{1}{T-t+\delta^{2}} \biggr)
\\
&&\quad\qquad{} + \biggl(\frac{-a'_{t}}{x-a_{t}+\delta}+\frac{b'_{t}}{b_{t}-x+\delta
} \biggr)
\\
&&\qquad \leq \cases{ \displaystyle-\frac{1}{4}\frac{1}{\delta^{2}\wedge
d_{\mathcal{D}} (t,x )^{2}}, &\quad
whenever $d_{\mathcal{D}} (t,x )\leq\eta$,
\vspace*{2pt}\cr
c_{1}, &\quad
otherwise.}
\end{eqnarray*}
Now set
\[
v^{2} (t,x )= \biggl(\frac{1}{\eta^{2}}+c_{1} \biggr)
\Bigl(\sup_{s\in(0,T )}a_{s}-x \Bigr) \Bigl( \Bigl(\inf
_{s\in(0,T )}b_{s} \Bigr)-x \Bigr).
\]
It follows that $v^{2}\geq0$ on $\mathcal{D}$ and
\[
\biggl(\partial_{t}+\frac{1}{2}\partial_{xx}
\biggr)v^{2}=- \biggl(\frac{1}{\eta^{2}}+c_{1} \biggr).
\]
Hence choosing
\[
v=v^{1}+v^{2}+3\llvert\log\delta\rrvert
\]
and putting the above together implies $ (\partial_{t}+\frac
{1}{2}\partial_{xx} )v\leq-\frac{c_{2}}{d_{\mathcal{D}}^{2}\wedge\delta^{2}}$
on $\mathcal{D}$. Since
\[
d_{\mathcal{D}}^{2} (s,y )\leq c_{3}d_{\mathcal{D}}^{2}
(t,x )
\]
for all $ (s,y )\in R_{t,x}^{\rho(t,x )}$, we
obtain from (\ref{eqPVIto}) that for all $ (t,x )\in\mathcal
{D}\setminus\mathcal{D}_{\delta}$,
\begin{eqnarray*}
(Pv-v ) (t,x ) & \leq& -\frac{c_{2}}{c_{3}d_{\mathcal{D}}^{2} (t,x
)}\mathbb{E} \bigl[
\tau^{t,x}-t \bigr]
\\
& =&  -\frac{c_{2}\rho^{2} (t,x )}{c_{3}d_{\mathcal{D}}^{2} (t,x )}
\\
&\leq& -\frac{1}{C}.
\end{eqnarray*}
Since in addition $v\geq0$ on $\overline{\mathcal{D}}$, it follows by
comparison with (\ref{eqPn}) that the expected number of steps satisfies
\begin{eqnarray*}
n (t,x ) & \leq& Cv (t,x )\leq C \bigl(1+\llvert\log\delta\rrvert\bigr).
\end{eqnarray*}\upqed
\end{pf}

%ex1 #&#
\begin{example}
\label{exrootMC}To give a numerical example, consider the function
\[
u (t,x )=4x^{4}+24 (1-t )x^{2}+12 (1-t )^{2}.
\]
It is a simple explicit solution of the unrestricted heat equation,
and by setting
\[
g (t,x )=u (t,x )
\]
on the parabolic boundary, it becomes the unique $C^{1,2}$ solution
of (\ref{eqpPDE}). In Figure~\ref{fig4} are the numerics for the choice
\[
T=1,\qquad
a_{t}=2-t,\qquad
b_{t}=0
\]
and $\rho(t,x )=\min(\frac{2-t-x}{\sqrt{2}},1-t,x )$
for $u (0,1 )=40$.
%
%f4 #&#
\begin{figure}[t]

\includegraphics{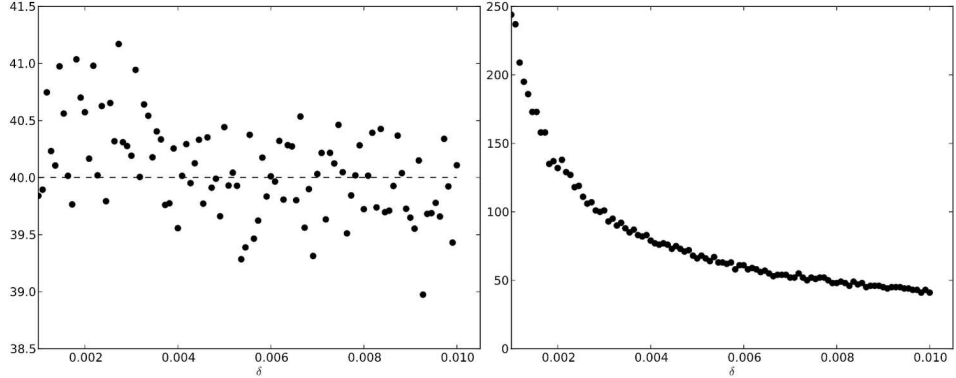}

\caption{The results for the Random walk over Root barriers applied to Example
\protect\ref{exrootMC}. The figure on the left shows the approximation to
$u (0,1 )=40$ and the right-hand figure the average number
of steps taken before leaving the domain, both as function of $\delta
\in\{ 0.001,\ldots,0.01 \}$. Each point
represents a run of the Monte Carlo scheme with 10,000 samples trajectories.}\label{fig4}
\end{figure}
\end{example}

%s6 #&#
\section{Conclusion and possible extensions}

We have presented a new characterisation of Root's solution of the
classic Skorokhod embedding problem (\ref{SEP}) by identifying it as
the unique solution
of an integral equation that has an intuitive interpretation and simple
derivation.
We then provided conditions on $\mu$ which imply geometric properties
about the shape of the barrier.
This in turn simplifies the integral equation for numerical purposes.
Finally, we have shown that the Root barrier can be used to yield a new
and very simple random walk over spheres algorithm.
It is natural to ask for several extensions:
\begin{itemize}
\item The proof of Theorem~\ref{teointegral} can be extended to other
processes than Brownian motion. While existence of the Root barrier
is known, the issue is to find explicit formulas for the expected
local time of this process to make this actually useful for numerics
(note that this is not needed for the PDE approach). Similarly,
Section~\ref{secintegralequatoin} applies (with minor modifications) to
the case of one-dimensional Brownian motion started with any probability
measure that is in convex order with $\mu$.
\item Not much is known about (\ref{SEP}) in multi dimensions.\footnote
{For the Root solution some existence results are known \cite
{MR0445600} but do not apply immediately;
for example, one-point sets are not regular anymore which leads to
issues about
randomised stopping times, etc.%
} However, for radially symmetric target measures (like the uniform
distribution on the unit ball) and multidimensional Brownian motion,
the question is equivalent to embedding into the Bessel process; hence
one can apply a simple modification of Theorem~\ref{teointegral}
in which the expected local time has still an explicit form. Unfortunately,
for the general multidimensional (or even non-Brownian) case, new ideas
are needed, and we hope to return to this and related Monte Carlo applications
in future work.
\item Section~\ref{secbarrierfunction} provides sufficient conditions
on $\mu$ such that its barrier function becomes monotone, and the
integral equation (\ref{eqintequation}) simplifies to a Volterra
equation of the first kind. Numerics for nonlinear integral equations
are a well-studied topic, and in principle one could hope to find fast
numerics for the integral equation (\ref{eqintequation}) such that
also for the general atom-free target measures equation (\ref{eqintequation})
becomes a competitor in numerics to the nonlinear PDE approach.
\end{itemize}

\section*{Acknowledgement}
Harald Oberhauser would like to thank the organisers and participants of the 6th
European Summer school in Financial Mathematics in Vienna 2013 for
helpful remarks.
%Harald Oberhauser and Paul Gassiat are grateful for partial support from the
%European Research Council, Grant agreement Nr.~258237.
%Harald Oberhauser is grateful for partial support from the European Research Council,
%Grant agreement Nr.~291244.

%\begin{appendix}
%\section{}
%\end{appendix}

%\section*{Acknowledgments}

%\begin{supplement}[id=suppA]
%\slink[doi]{10.1214/00-AAPXXXXSUPP} %[doi,text={...}] - jei reikia
%suskaldyti doi
%\sfilename{aapXXXX\_supp.pdf}
%\end{supplement}

% imsref loaded by linak, 2014-07-04 13:20:40
%
% imsref loaded by linak, 2014-08-14 15:15:25

\printaddresses
\end{document}